\documentclass[12pt]{amsart}
\input epsf
\makeatletter
\let\atopwithdelims\@@atopwithdelims
\let\over\@@over
\newtheorem{theorem}{Theorem}[section]
\newtheorem{proposition}[theorem]{Proposition}
\newtheorem{corollary}[theorem]{Corollary}

\newcounter{enumit}
\newenvironment{enumit}{\begin{list}{({\it\roman{enumit}})}{
                        \usecounter{enumit}}
                      }{\end{list}}

\newcommand{\toto}{\longrightarrow}

\newcommand{\for}{{\rm for}}
\newcommand{\co}{{\bf c}}
\begin{document}
\title[Piles of Cubes and Monotone Path Polytopes]{Piles of Cubes, 
Monotone Path Polytopes and Hyperplane Arrangements}
\author{Christos~A.~Athanasiadis}
\address{\hskip-\parindent Christos~A.~Athanasiadis\\
Mathematical Sciences Research Institute\\
1000 Centennial Drive\\
Berkeley, CA 94720}
\email{athana@msri.org}
\thanks{Supported by a postdoctoral fellowship from the Mathematical 
Sciences Research Institute, Berkeley, California. Research at MSRI 
is supported in part by NSF grant DMS-9022140.}
\begin{abstract}
Monotone path polytopes arise as a special case of the construction of fiber polytopes, introduced by Billera and Sturmfels. A simple example is provided by the permutahedron, which is a monotone path polytope of the standard unit cube. The permutahedron is the zonotope polar to the braid arrangement. We show how the zonotopes polar to the cones of certain deformations of the braid arrangement can be realized as monotone path polytopes. The construction is an extension of that of the permutahedron and yields interesting connections between enumerative combinatorics of hyperplane arrangements and geometry of monotone path polytopes.
\end{abstract}

\maketitle

\section{Introduction}

Fiber polytopes were introduced by Billera and Sturmfels \cite{BS1} and were further studied in \cite{BS2}. The fiber polytope $\Sigma(P,Q)$ is a polytope naturally associated to any projection of polytopes $\pi : P \toto Q$. It is defined as
\begin{equation}
\Sigma(P,Q) = \left\{ \frac{1}{{\rm vol}(Q)} \int_Q \gamma(x) \, {\rm d}x \ | \ \textrm{$\gamma$ is a section of $\pi$} \right\}. 
\label{sec}
\end{equation}
A concise introduction to the theory of fiber polytopes can be found in Chapter 9 of Ziegler's book \cite{Zi}, on which we rely for general background on polytopes. 

We are only interested in the special case in which $Q$ is one dimensional, i.e. a line segment. The polytope $\Sigma(P,Q)$ is then called the \textit{monotone path polytope} of $P$ and $\pi$ \cite[\S 5]{BS1} \cite[\S 9.2]{Zi} and is denoted by $\Pi(P, \pi)$. Its dimension is one less than that of $P$. The sections of $\pi$ are the paths in $P$ which are \textit{monotone}, that is strictly increasing with respect to the function $\pi$, and \textit{maximal}, i.e. they project to $Q$. The vertices of $\Pi(P, \pi)$ correspond to certain maximal monotone edge paths on $P$, called $\pi$-\textit{coherent}. 

\begin{figure}[htpb]
\center{\mbox{\epsfbox{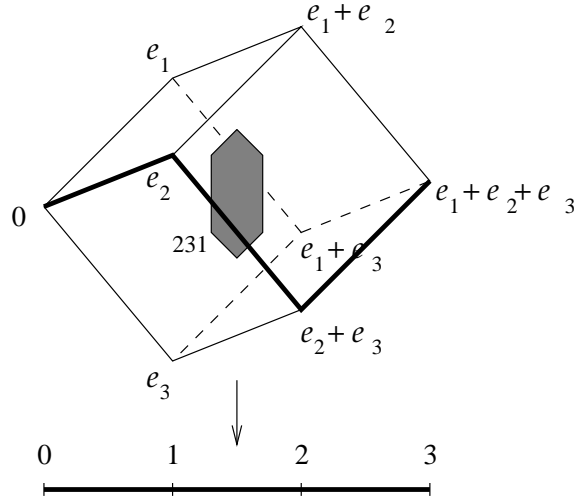}}} 
\caption{The projection $\pi$ and the permutahedron}
\label{pi}
\end{figure}

\vspace{0.1 in}
A nice example of a monotone path polytope is the permutahedron ${\Pi}_{d-1} \subseteq {\mathbb R}^d$. The permutahedron ${\Pi}_{d-1}$ is a classical geometric object. It is the convex hull of all vectors in ${\mathbb R}^d$ obtained by permuting the coordinates of the vector $(1, 2,\ldots,d)$. It was constructed as a monotone path polytope of the standard unit cube $C_d = [0, 1]^d$ by Billera and Sturmfels \cite[Example 5.4]{BS1} (see also \cite[Thm.\ 4.3]{BS2}). A detailed description of the construction can be found in \cite[Example 9.8]{Zi}. The projection used is $\pi : C_d \toto [0,d]$ with
\begin{equation}
\pi(x) = x_1 + x_2 + \cdots + x_d. 
\label{sum}
\end{equation}
The monotone edge paths from the origin to the vertex $(1,1,\ldots,1)$ are all $\pi$-coherent. They correspond naturally to permutations of a $d$-element set and give rise to the vertices of ${\Pi}_{d-1}$. This is illustrated in Figure \ref{pi} for $d=3$.

\begin{figure}[htpb]
\center{\mbox{\epsfbox{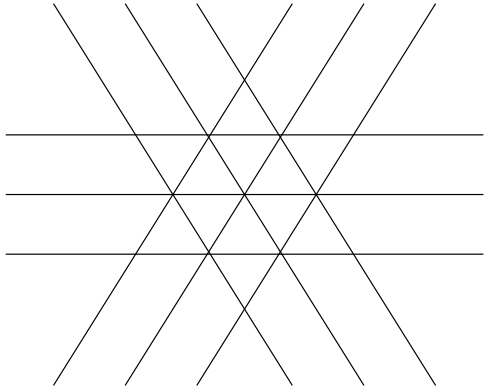}}} 
\caption{The Catalan arrangement for $d=3$}
\label{C3}
\end{figure}

An important property of ${\Pi}_{d-1}$ is that it is a zonotope. The associated hyperplane arrangement \cite[\S 7.3]{Zi} is the \textit{braid arrangement} ${\mathcal A}_d$. It consists of the hyperplanes in ${\mathbb R}^d$ of the form $x_i = x_j$ for $1 \leq i < j \leq d$, i.e. the reflecting hyperplanes of the Coxeter group of type $A_{d-1}$. Our motivation comes from recent work on \textit{deformations} of ${\mathcal A}_d$, initiated by Stanley \cite{St1}. A deformation of ${\mathcal A}_d$ is an affine arrangement which has each of its hyperplanes parallel to some hyperplane $x_i = x_j$ of ${\mathcal A}_d$. We will be interested in the deformations of the form 
\begin{equation}
x_i - x_j = - \lambda_i + 1,\ldots,-1, 0, 1,\ldots,\lambda_j - 1 \ \ 
\for \ \ 1 \leq i < j \leq d,
\label{arr}
\end{equation}
where $\lambda = (\lambda_1, \lambda_2,\ldots,\lambda_d)$ is a composition of a positive integer $n$. We denote the arrangement (\ref{arr}) by ${\mathcal A}_d (\lambda)$. For $\lambda_1 = \cdots = \lambda_d = a+1$ it specializes to the \textit{extended Catalan arrangement} (see for example \cite[\S 2]{St1}) and is denoted by ${\mathcal A}_d ^{[0,a]}$. This reduces to the braid arrangement ${\mathcal A}_d$ for $a=0$ and is simply called a \textit{Catalan arrangement} for $a=1$. Figure \ref{C3} shows the Catalan arrangement for $d=3$, intersected with the hyperplane $x_1 + x_2 + x_3 = 0$.

\vspace{0.1 in}
The combinatorics of deformations of ${\mathcal A}_d$ is very rich and has a growing literature \cite{St1, St2, Ath3, PS}. The zonotopes associated to the corresponding homogenized linear arrangements, or cones, haven't been studied much as polytopes. They are generalizations of the permutahedron that carry the same interesting combinatorial structure as the arrangements. Our objective is to construct the zonotopes associated to the cones of the arrangements ${\mathcal A}_d (\lambda)$, with their Minkowski summands suitably rescaled, as monotone path polytopes. 

The construction generalizes that of the permutahedron. In particular, the cube $C_d$ is replaced by Ziegler's \textit{lifted pile of cubes} $\widetilde{\mathcal P}_{d+1} (\lambda)$ \cite[\S 5.1]{Zi}, which we define formally in the next section. The projection $\pi$ is given by the formula (\ref{sum}). The corresponding monotone path polytope $\Pi_d (\lambda)$ is a zonotope, as described above, with its upper part truncated. The maximal monotone edge paths in the lower faces of $\widetilde{\mathcal{P}}_{d+1} (\lambda)$ are in bijection with the lattice paths in ${\mathbb R}^d$ from the origin to the point $\lambda$, having unit coordinate steps. In general, most of these paths are \textit{not} coherent. The coherent ones correspond to the lower vertices of the zonotope and hence to the regions of ${\mathcal A}_d (\lambda)$. This will enable us to enumerate them (see Corollary 4.3) by an easy application of the ``finite field method'' of \cite{Ath1, Ath2}. It follows that their number is small compared to the total number of lower monotone edge paths.

\vspace{0.1 in}
The paper is organized as follows: In Section 2 we give the necessary definitions and state our main results about $\Pi_d (\lambda)$. Section 3 contains a proof of the main theorem. In Section 4 we draw all enumerative consequences of the main theorem by analyzing the combinatorics of ${\mathcal A}_d (\lambda)$. We also characterize the coherent paths in the lower faces of $\widetilde{\mathcal P}_{d+1} (\lambda)$ directly from the definition. We close with some remarks in Section 5.

\section{Basic definitions and results}

We begin with some notation and terminology. We denote by $e_1, e_2,\ldots,e_d$ the unit coordinate vectors in ${\mathbb R}^d$. We write $[a, b]$ for the line segment joining two points $a, b \in {\mathbb R}^d$. The \textit{canonical projection} ${\mathbb R}^{d+1} \toto {\mathbb R}^d$ is the map which simply forgets the last coordinate. Lastly, we refer to the union of the lower faces \cite[\S 5.1]{Zi}, or upper faces, of a polytope $P$ as its \textit{lower part}, or \textit{upper part} respectively.

\begin{figure}[htpb]
\center{\mbox{\epsfbox{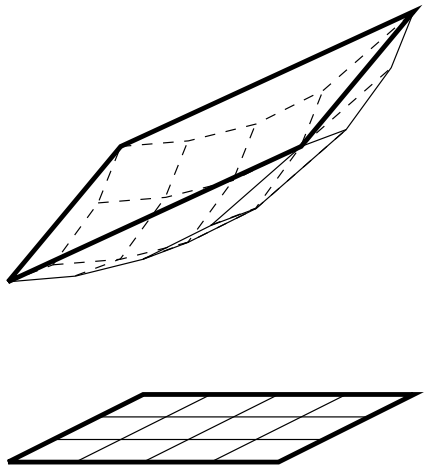}}} 
\caption{${\mathcal P}_2 (4,3)$ and $\widetilde{\mathcal P}_3 (4,3)$}
\label{pile}
\end{figure}

We first give the definition of the polytope $\widetilde{\mathcal P}_{d+1} (\lambda)$, as promised in the introduction. Let $\lambda = (\lambda_1, \lambda_2,\ldots,\lambda_d)$ be a composition of a positive integer $n$. In other words, the $\lambda_i$ are positive integers which sum to $n$. The \textit{pile of cubes} ${\mathcal P}_d (\lambda)$, corresponding to $\lambda$, is the polytopal complex formed by all unit cubes with integer vertices in the \textit{d-box}  
\[ B(\lambda) = \left\{ (x_1,\ldots,x_d) \in {\mathbb R}^d \ | \ 0 \leq x_i 
\leq \lambda_i \ {\rm for} \ 1 \leq i \leq d \right\}. \]  
In particular, the vertex set of ${\mathcal P}_d (\lambda)$ is
\[ {\rm vert} \left({\mathcal P}_d (\lambda) \right) = B(\lambda) \cap {\mathbb Z}^d. \]
Let $f : {\mathbb R}^d \toto {\mathbb R}$ be of the form 
\begin{equation}
f(x) = f_1(x_1) + \cdots + f_d(x_d)
\label{sep}
\end{equation}
where the $f_i$'s are strictly convex functions of one variable. For the purpose of our results we will use the canonical choice
\[ f(x) = x_1^2 + x_2^2 + \cdots + x_d^2. \]
We take this as the definition of $f$ from now and on. The \textit{lifted pile of cubes} corresponding to $\lambda$ is the polytope in ${\mathbb R}^{d+1}$
\[ \widetilde{\mathcal P}_{d+1} (\lambda) = {\rm conv} \left\{ (x, f(x)) \ | \ x \in {\rm vert} \left({\mathcal P}_d (\lambda) \right) \right\}. \]
The pile of cubes ${\mathcal P}_d (\lambda)$ and its lift $\widetilde{\mathcal P}_{d+1} (\lambda)$ are discussed in \cite[Example 5.4]{Zi}. For suitable $\lambda$ they provide examples of polytopal complexes, respectively polytopes, which are not extendably shellable \cite[Ch.\ 8]{Zi}. Figure \ref{pile} shows ${\mathcal P}_2 (4,3)$ and its lift $\widetilde{\mathcal P}_3 (4,3)$. 

The lower faces of $\widetilde{\mathcal P}_{d+1} (\lambda)$ form a polytopal complex combinatorially equivalent to ${\mathcal P}_d (\lambda)$. Because of the form (\ref{sep}), these faces are parallelepipeds, as Figure \ref{pile} suggests. By construction, they map to the corresponding faces of ${\mathcal P}_d (\lambda)$ under the canonical projection ${\mathbb R}^{d+1} \toto {\mathbb R}^d$. The unique upper facet, defined by the hyperplane
\begin{equation}
\lambda_1 x_1 + \lambda_2 x_2 + \cdots + \lambda_d x_d = x_{d+1},
\label{hyp}
\end{equation}
is a parallelepiped which projects to $B(\lambda)$.

\vspace{0.1 in}
Let $\Pi_d (\lambda)$ be the monotone path polytope of the projection 
\[ \pi : \widetilde{\mathcal P}_{d+1} (\lambda) \toto [0, n] \]
given by the formula (\ref{sum}). Thus $\pi(x)$ is the sum of the first $d$ coordinates of $x \in {\mathbb R}^{d+1}$. We have $\dim \Pi_d (\lambda) = \dim \widetilde{\mathcal P}_{d+1} (\lambda) - 1= d$, the number of parts of $\lambda$ except if $\lambda_i = 1$ for all $i$. In this case $\widetilde{\mathcal P}_{d+1} (\lambda)$ is affinely isomorphic to the $d$-cube $C_d$, so $\Pi_d (\lambda)$ is affinely isomorphic to the permutahedron $\Pi_{d-1}$ and $\dim \Pi_d (\lambda) = d-1$.

\vspace{0.1 in}
A \textit{truncation} of a polytope $P$ in ${\mathbb R}^r$ is its intersection with a closed halfspace $a \cdot x \leq 0$, where $a \in {\mathbb R}^r$. We call such a truncation an \textit{upper truncation} if it includes the lower part of $P$ but does not intersect the upper part in the region $a \cdot x < 0$. In particular, an upper truncation has a unique upper facet. It rarely happens that a polytope $P$ has an upper truncation. However, if $P$ is a zonotope with a sufficiently large Minkowski summand in the direction of the last coordinate, as in the situation of the next theorem, then clearly, upper truncations exist. 
\begin{theorem}
Let $Z_d (\lambda)$ be the zonotope
\[ \left[ 0, s \, e_{d+1} \right] + \sum_{1 \leq i < j \leq d} 
\sum_{k=0}^{\lambda_i - 1} \sum_{l=0}^{\lambda_j - 1}
\left[ 0, e_j - e_i + (2l-2k) \, e_{d+1} \right], \]
where $s$ is a sufficiently large positive number. The polytope $\Pi_d (\lambda)$ is an upper truncation of a translate of $\frac{1}{n} Z_d (\lambda)$, cut by the hyperplane (\ref{hyp}).
\end{theorem}

\begin{figure}[htpb]
\center{\mbox{\epsfbox{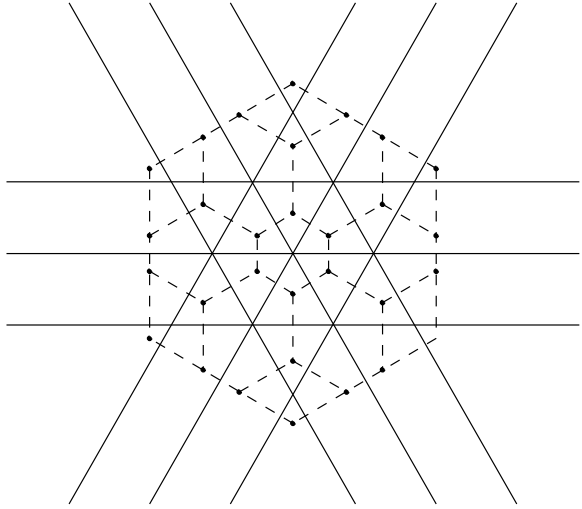}}} 
\caption{The polar complex}
\label{pol}
\end{figure}

Note that the hyperplane arrangement associated to $Z_d (\lambda)$ is projectively equivalent to the \textit{cone} of ${\mathcal A}_d (\lambda)$, denoted $\co {\mathcal A}_d (\lambda)$. This is the linear arrangement in ${\mathbb R}^{d+1}$ obtained by homogenizing each hyperplane $x_i - x_j = s$ of ${\mathcal A}_d (\lambda)$ to $x_i - x_j = s x_{d+1}$ and adding the hyperplane $x_{d+1} = 0$. Hence $Z_d (\lambda)$ is combinatorially equivalent to the zonotope $Z[\co {\mathcal A}_d (\lambda)]$, associated to this cone, which we call the \textit{polar zonotope}. The faces of $Z[\co {\mathcal A}_d (\lambda)]$ are in inclusion-reversing bijection with the faces of $\co {\mathcal A}_d (\lambda)$ \cite[Cor.\ 7.17]{Zi}. Its lower faces define a zonotopal complex $Z[{\mathcal A}_d (\lambda)]$ which we call the \textit{polar zonotopal complex}, or simply the \textit{polar complex}. The faces of $Z[{\mathcal A}_d (\lambda)]$ are in inclusion-reversing bijection with those of ${\mathcal A}_d (\lambda)$. Figure \ref{pol} shows the complex polar to the Catalan arrangement of Figure \ref{C3}, projected canonically on ${\mathbb R}^2$.

\vspace{0.1 in}
We now focus on the the lower part, which is the interesting part of $\Pi_d (\lambda)$. Theorem 2.1 implies the following corollary. 
\begin{corollary} 
The complex of all lower faces of $\Pi_d (\lambda)$ is combinatorially equivalent to the polar complex $Z[{\mathcal A}_d (\lambda)]$ of the arrangement ${\mathcal A}_d (\lambda)$. In particular, the lower vertices of $\Pi_d (\lambda)$ are in bijection with the regions of ${\mathcal A}_d (\lambda)$.
\qed
\end{corollary}
Hence, for $\lambda = (2,2,2)$, $\Pi_d (\lambda)$ looks combinatorially like the complex of Figure \ref{pol}, from a point far below on the $x_3$ axis. 

\vspace{0.1 in}
It follows from Corollary 2.2 that the $\pi$-coherent monotone edge paths in the lower part of $\widetilde{\mathcal P}_{d+1} (\lambda)$ are in bijection with the regions of ${\mathcal A}_d (\lambda)$. To characterize them, we introduce some more terminology. We call the maximal monotone edge paths (with respect to $\pi$) in ${\mathcal P}_d (\lambda)$ simply \textit{monotone $\lambda$-paths}. They are the lattice paths in ${\mathbb R}^d$ from the origin to the point $\lambda$ having $n$ unit coordinate steps. They correspond naturally to permutations of the multiset $M_{\lambda}$, which contains $i$ with multiplicity $\lambda_i$ for $1 \leq i \leq d$. We call these permutations $\lambda$-\textit{permutations} and denote by $p^w$ the monotone $\lambda$-path which corresponds to $w$. For example, if $\lambda = (2, 2, 2)$ and $w = 1 2 1 3 2 3$ then the path $p^w$ has successive steps $e_1, e_2, e_1, e_3, e_2, e_3$. The number of $\lambda$-permutations is the multinomial coefficient $n \choose \lambda_1,\ldots,\lambda_d$.

The maximal monotone edge paths in the lower part of $\widetilde{\mathcal P}_{d+1} (\lambda)$ are in bijection with the monotone $\lambda$-paths, via the canonical projection ${\mathbb R}^{d+1} \toto {\mathbb R}^d$. We refer to these paths in $\widetilde{\mathcal P}_{d+1} (\lambda)$ as \textit{lifted monotone $\lambda$-paths} and use the notation $\gamma^w$. 

\vspace{0.1 in}
A $\lambda$-permutation $w = w_1 w_2 \cdots w_n$ is \textit{nesting} if there exist indices $1 \leq i < j < k < l \leq n$ such that $w_i = w_l$, $w_j = w_k$ and $w_m \neq w_i$ for all $m$ with $j \leq m \leq k$. Otherwise $w$ is \textit{non-nesting}. Equivalently, $w$ is non-nesting if the following linear system in the variables $y_1, y_2,\ldots,y_n$ is feasible: 
\begin{enumit}
\item
$y_1 < y_2 < \cdots < y_n$,
\item
$y_k = y_j + 1$ \ if \ $1 \leq j < k \leq n$, $w_k = w_j$ and $w_m \neq w_j$ 
for $j < m < k$.
\end{enumit}
For example, $1 2 1 2 1$ is non-nesting but $1 2 2 1 1$ is not.
\begin{theorem} 
The lifted monotone $\lambda$-path $\gamma^w$ is $\pi$-coherent if and only if $w$ is non-nesting. The number of such paths, or equivalently the number of lower vertices of $\Pi_d (\lambda)$, or regions of ${\mathcal A}_d (\lambda)$ is $\frac{n!}{(n-d+1)!}$. In particular, this quantity depends only on $n$ and $d$.
\end{theorem}

It follows that the fraction of monotone $\lambda$-paths with $\pi$-coherent lifts is
\[ \frac{\lambda_1 ! \lambda_2 ! \cdots \lambda_d !}{(n-d+1) !}. \]
This result gives concrete examples of projections for which the fraction of coherent monotone edge paths is explicitly shown to be small. For instance, in the case $\lambda_1 = \cdots = \lambda_d = 2$ this quantity becomes $2^d/(d+1)!$. The only other fiber polytope situation we are aware of, in which it has been shown that the number of coherent subdivisions is asymptotically negligible compared to the number of all subdivisions, appears in \cite{dHSS} in the context of triangulations of cyclic $d$-polytopes with $d+4$ vertices.

\section{Proof of the main theorem}

We now prove Theorem 2.1. 

\vspace{0.1 in}
\noindent
\textit{Proof of Theorem 2.1}: 
We use directly the defining formula (\ref{sec}). Recall that the vertices of $\Pi_d (\lambda)$ correspond to certain maximal monotone edge paths $\gamma$ in $\widetilde{\mathcal P}_{d+1} (\lambda)$. Since $\pi (x)$ depends only on the first $d$ coordinates of $x$, we only need to consider such paths either in the lower or upper part of $\widetilde{\mathcal P}_{d+1} (\lambda)$ and take the convex hull of the right hand side of (\ref{sec}). Let
\[ I_{\gamma} := \frac{1}{\textrm{vol} (Q)} \int_{[0, n]} \gamma(x) \, 
{\rm d} x = \frac{1}{n} \, \sum_{k=1}^n \frac{\gamma(k-1) + \gamma(k)}{2}.
\]
We first consider the lower part of $\widetilde{\mathcal P}_{d+1} (\lambda)$,
so the relevant paths are the lifted monotone $\lambda$-paths $\gamma^w$. 
Instead of calculating $I_{\gamma^w}$ explicitly, we note that any $\lambda$-permutation $w$ can be obtained from $w_0 = 1 \cdots 1 \cdots d \cdots d$ by successively swapping adjacent entries $i j$ with $i < j$. Let $w = w_1 w_2 \cdots w_n$, where $w_k = i < j = w_{k+1}$ and let $w^{\prime}$ be obtained from $w$ by swapping $w_k$ and $w_{k+1}$. If $v := \gamma^w (k-1) = \gamma^{w^{\prime}} (k-1) = (m_1,m_2,\ldots,m_d)$, then
\[ I_{\gamma^{w^{\prime}}} - I_{\gamma^w} = {\textstyle \frac{1}{n}} \left( 
\gamma^{w^{\prime}} (k) - \gamma^w (k) \right) = \]
\[ {\textstyle \frac{1}{n}} \, \left( v + e_j + f(v+e_j) \, 
e_{d+1} - v - e_i - f(v+e_i) \, e_{d+1} \right) = \]
\[ {\textstyle \frac{1}{n}} \, \left( e_j - e_i + (2m_j - 2m_i) \, 
e_{d+1} \right). \]

\vspace{0.1 in}
\noindent
It follows easily that the points $I_{\gamma^w}$ lie on $\frac{1}{n} Z_d (\lambda) + I_{\gamma^{w_0}}$ and include all its lower vertices.

Finally suppose that $\gamma$ lies in the upper facet of $\widetilde{\mathcal P}_{d+1} (\lambda)$, which is defined by the hyperplane (\ref{hyp}). This upper facet projects to the $d$-box $B(\lambda)$ under the canonical projection ${\mathbb R}^{d+1} \toto {\mathbb R}^d$. Each of the monotone edge paths $\gamma$ projects to a monotone $\lambda$-path with corresponding $\lambda$-permutation of the form $w = \sigma_1 \cdots \sigma_1 \cdots \sigma_d \cdots \sigma_d$, where $\sigma_1 \cdots \sigma_d$ is a permutation of $\{1,2,\ldots,d\}$. Hence $I_{\gamma} = I_{\gamma^w} + t e_{d+1}$ for some $t \geq 0$. Since $\gamma$ lies on (\ref{hyp}), so does $I_{\gamma}$. It follows that the upper part of $\Pi_d (\lambda)$ also consists of a single facet, defined by (\ref{hyp}). It projects to the same zonotope on ${\mathbb R}^d$, which is combinatorially equivalent to the permutahedron $\Pi_{d-1}$, as the lower part. This completes the proof.
\qed

\section{The coherent paths and enumeration}

In this section we characterize and enumerate the $\pi$-coherent lifted monotone $\lambda$-paths and thus prove Theorem 2.3. We proceed directly from the definition of a coherent path and extend our argument to get another proof of Corollary 2.2. For the enumeration we use the finite field method of \cite{Ath1, Ath2}.

\vspace{0.1 in}
We first recall what it means for a monotone edge path to be coherent. Let $P$ be a polytope in ${\mathbb R}^r$, $\pi : {\mathbb R}^r \toto {\mathbb R}$ an affine function, $Q = \pi(P) \subseteq {\mathbb R}$ and $p$ a maximal monotone edge path in $P$. The path $p$ is said to be $\pi$-\textit{coherent} if there exists a \textit{generic} linear functional $c \in ({\mathbb R}^r)^{\ast}$ such that $\pi^c : {\mathbb R}^r \toto {\mathbb R}^2$ with
\[ \pi^c (x) = \left(\pi(x), c(x)\right) \]
maps $p$ to the path of lower edges of the polygon $Q^c := \pi^c(P)$. More generally \cite[Definition 9.2]{Zi}, a set of faces $\mathcal F$ of $P$ is said to form a $\pi$-\textit{coherent subdivision} of $Q$ if, for some $c \in ({\mathbb R}^r)^{\ast}$, $\mathcal F$ is the set ${\mathcal F}^c$ of faces $\left( \pi^c \right)^{-1} (G)$ of $P$, where $G$ runs through the lower faces of $Q^c$. Note that ${\mathcal F}^c$ induces the subdivision $\{ \pi(F) : F \in {\mathcal F}^c \}$ of $Q$. This is also the subdivision obtained by projecting canonically the lower edges of $Q^c$ on $Q$. The set ${\mathcal F}^c$ consists of faces of $P$, rather than $Q$.

\vspace{0.1 in}
For the rest of this section let $r=d+1$, $P = \widetilde{\mathcal P}_{d+1} (\lambda)$, $Q = [0, n]$ and $\pi$ as in Section 2. 
\begin{proposition} 
The lifted monotone $\lambda$-path $\gamma^w$ is $\pi$-coherent if and only if $w$ is non-nesting. The non-nesting $\lambda$-permutations are in bijection with the regions of the arrangement ${\mathcal A}_d (\lambda)$. 
\end{proposition}
\begin{proof}
Let $c \in ({\mathbb R}^{d+1})^{\ast}$ be generic and given by
\begin{equation}
c(x) = a_1 x_1 + a_2 x_2 + \cdots + a_{d+1} x_{d+1}.
\label{c}
\end{equation}
We need to describe the $\pi$-coherent path defined by $c$, directly in terms of $c$. Since we want this path to lie in the lower part of $\widetilde{\mathcal P}_{d+1} (\lambda)$, we can assume that $a_{d+1}$ is positive, say $a_{d+1} = 1/2$ by rescaling. 

Let $w = w_1 w_2 \cdots w_n$ be a $\lambda$-permutation and $p^w, \gamma^w$ as in Section 2, with $p^w = (v_0, v_1,\ldots,v_n)$ and $\gamma^w = (\tilde{v}_0, \tilde{v}_1,\ldots,\tilde{v}_n)$. Thus $v_k = v_{k-1}+e_i$ if $w_k = i$ and $\tilde{v}_k = (v_k, f(v_k))$ for all $k$. We denote by $\Delta^c (w)$ the sequence of length $n$ whose $k$th term records the difference
\[ c(\tilde{v}_k) - c(\tilde{v}_{k-1}) \]
of $c$ at the $k$th edge of $\gamma^w$. Let $v_{k-1} = (m_1,\ldots,m_d)$ and $w_k = i$. Then
\[ c(\tilde{v}_k)-c(\tilde{v}_{k-1}) = a_i + a_{d+1} \, (f(v_{k-1} + e_i) - 
f(v_{k-1})) = \]
\[ a_i + (1/2) \, (2m_i + 1) = a^{\prime}_i + m_i, \]

\noindent
where $a^{\prime}_i = a_i + 1/2$ for $1 \leq i \leq d$ and the possible values of $m_i$ are the integers satisfying $0 \leq m_i \leq \lambda_i - 1$. It follows that $\Delta^c (w)$ is the permutation of the numbers
\begin{equation}
a^{\prime}_i + m, \ \ 0 \leq m \leq \lambda_i - 1, \ 1 \leq i \leq d 
\label{nu}
\end{equation}
which places $a^{\prime}_i, a^{\prime}_i + 1,\ldots,a^{\prime}_i + \lambda_i - 1$, in this order, in the $\lambda_i$ positions occupied by $i$ in $w$. The numbers (\ref{nu}) are distinct because $c$ is assumed to be generic.

Let $w^c$ be the $\lambda$-permutation which corresponds to the $\pi$-coherent path defined by $c$. If $\pi^c$ maps $\gamma^w$ to the lower edge path of the polygon $Q^c$, then $\Delta^c (w)$ is strictly increasing, by the convexity of this path. It follows that $w^c$ is the unique $\lambda$-permutation $w$ for which the permutation $\Delta^c (w)$ of the numbers (\ref{nu}) is strictly increasing. Let $c^{\prime} = (a_1^{\prime}, a_2^{\prime},\ldots,a_d^{\prime}) \in {\mathbb R}^d$. For a fixed $\lambda$-permutation $w$, if nonempty, the set of all $c^{\prime}$ such that $w^c = w$ forms exactly a region of ${\mathcal A}_d (\lambda)$. Indeed, deciding which of the two numbers $a_i^{\prime} + m$ and $a_j^{\prime} + l$ in (\ref{nu}) is smaller amounts to choosing a side of the hyperplane $x_i - x_j = l - m$ of ${\mathcal A}_d (\lambda)$. Finally, $w=w^c$ for some $c$ if and only if the condition that $\Delta^c (w)$ is strictly increasing does not impose inequalities of the form
\[ d_j < d_i < d_i + 1 < d_j + 1. \]
By definition, this means that $w$ is non-nesting.
\end{proof}

\vspace{0.1 in}
We now compute the number of regions of ${\mathcal A}_d (\lambda)$. We assume familiarity with the \textit{characteristic polynomial} $\chi ({\mathcal A}, q)$ of ${\mathcal A}$ \cite[\S 2.3]{OT}, a fundamental combinatorial invariant of ${\mathcal A}$, and use Zaslavsky's theorem \cite{Za} for the number of regions of ${\mathcal A}$. The characteristic polynomial of ${\mathcal A}_d (\lambda)$ can be computed by an easy application of the finite field method of \cite{Ath1} \cite[Part II]{Ath2}. This method reduces the computation to a simple counting problem in a vector space over a finite field. The argument in the following proposition is similar to the one given in \cite[Thm.\ 5.1]{Ath1} for the special case of the extended Catalan arrangements.  
\begin{proposition} 
We have
\[ \chi ({\mathcal A}_d (\lambda), q) = q \prod_{j=n-d+1}^{n-1} (q-j). \]
\end{proposition}
\begin{proof} 
Let $q$ be a large prime number and let ${\mathbb F}_q$ denote the finite field of integers mod $q$. Theorem 2.2 in \cite{Ath1} (also \cite[Thm.\ 5.2.1]{Ath2}) implies that $\chi ({\mathcal A}_d (\lambda), q)$ counts the number of $d$-tuples $(x_1, x_2,\ldots,x_d) \in {\mathbb F}_q ^d$ satisfying
\[ x_i - x_j \neq - \lambda_i + 1,\ldots,-1, 0, 1,\ldots,\lambda_j - 1 \ \ 
\for \ \ 1 \leq i < j \leq d. \]
Equivalently, we want to choose the $x_i$ so that the classes mod $q$
\begin{equation}
x_i, x_i + 1,\ldots,x_i + \lambda_i - 1 \ \ \for \ \ 1 \leq i \leq d
\label{str}
\end{equation}
are distinct. To count these $d$-tuples we first cyclically permute the $d$ strings (\ref{str}) in $(d-1)!$ ways. Then we distribute $q-n$ indistinguishable boxes in the $d$ spaces between successive strings in $q-n+d-1 \choose d-1$ ways. The $q-n$ boxes distributed stand for the $q-n$ classes in ${\mathbb F}_q$ not of the form (\ref{str}). Finally we assign to $x_1$ a specific value in ${\mathbb F}_q$ in $q$ ways, say $x_1 = 0$. The other classes of (\ref{str}) $x_1 + 1,\ldots$ and the boxes are naturally assigned the values $1,\ldots,q - 1 \in {\mathbb F}_q$, according to their cyclic arrangement. The product 
\[ q \, (d-1)! \, {q-n+d-1 \choose d-1} \] 
is the expression for $\chi ({\mathcal A}_d (\lambda), q)$ we have claimed.
\end{proof}

\vspace{0.1 in}
Zaslavsky's theorem \cite{Za} expresses the number of regions of an arrangement ${\mathcal A}$ in ${\mathbb R}^d$ as $(-1)^d \chi ({\mathcal A}, -1)$ and yields the following corollary.
\begin{corollary} 
The number of regions of ${\mathcal A}_d (\lambda)$ is
\[ \frac{n!}{(n-d+1)!}. \]
In particular, this number depends only on $n$ and $d$, the sum and number of parts of $\lambda$ respectively.
\qed
\end{corollary}

Figure \ref{3,2,1} shows the arrangement ${\mathcal A}_d (\lambda)$ for $\lambda = (3, 2, 1)$, intersected with the hyperplane $x_1 + x_2 + x_3 = 0$. Note that it has $30$ regions, as many as the Catalan arrangement of Figure \ref{pol}, which corresponds to $\lambda = (2, 2, 2)$. Corollary 4.3 specializes to the formula given in \cite[\S 2]{St1} for the number of regions of the extended Catalan arrangement ${\mathcal A}_d ^{[0,a]}$, since in this case $n = (a+1)d$. Proposition 4.1 and Corollary 4.3 imply Theorem 2.3.

\begin{figure}[htpb]
\center{\mbox{\epsfbox{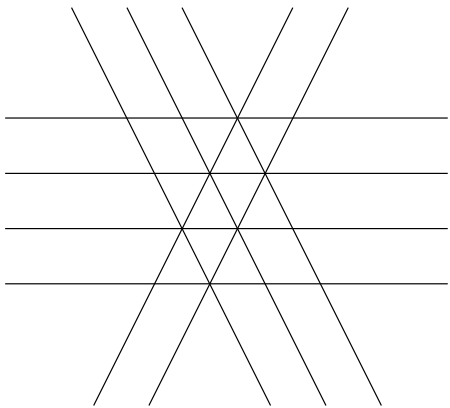}}} 
\caption{The arrangement ${\mathcal A}_3 (3, 2, 1)$}
\label{3,2,1}
\end{figure}

\vspace{0.1 in}
The reasoning in the proof of Proposition 4.1 can be extended to determine the face lattice of $\Pi_d (\lambda)$. Below we give an alternative proof of Corollary 2.2, obtained in this way. By the general theory of fiber polytopes \cite[Thm. 2.1]{BS1} \cite[Thm.\ 9.6]{Zi}, the face lattice of $\Pi_d (\lambda)$ is isomorphic to the poset of $\pi$-coherent subdivisions of $Q = [0, n]$, with a minimal element adjoined. In this poset we have ${\mathcal F}_1 \leq {\mathcal F}_2$ if the union of the faces in ${\mathcal F}_1$ is contained in the union of the faces in ${\mathcal F}_2$.

\vspace{0.1 in}
\noindent
\textit{Alternative proof of Corollary 2.2}: To get the lower faces of $\Pi_d (\lambda)$ we consider linear functionals $c \in ({\mathbb R}^{d+1})^{\ast}$ as in (\ref{c}), with positive $a_{d+1}$, but not necessarily generic. Let $P^c$ denote the set of monotone $\lambda$-paths whose lifts lie in ${\mathcal F}^c$, so that $P^c$ is a singleton if $c$ is generic. Note that ${\mathcal F}^{c_1} \leq {\mathcal F}^{c_2}$ in the poset of $\pi$-coherent subdivisions if and only if $P^{c_1} \subseteq P^{c_2}$. Using the notation in the proof of Proposition 4.1, we have $p^w \in P^c$ if and only if $\Delta^c (w)$ is \textit{weakly} increasing. It follows that $P^{c_1} \subseteq P^{c_2}$ if and only if the faces $F_1, F_2$ of ${\mathcal A}_d$ which contain $c^{\prime}_1$ and $c^{\prime}_2$ respectively satisfy $F_1 \geq F_2$ in the face poset of ${\mathcal A}_d$.
\qed

\vspace{0.1 in}
A precise description of the $\pi$-coherent lower subdivisions is given in Remark 2 of the next section. The computations of general face numbers of the extended Catalan arrangements \cite[Ch.\ 8]{Ath2} can be translated in terms of the lower part of $\Pi_d (\lambda)$ via Corollary 2.2. The following corollary, for example, follows from \cite[Cor.\ 8.3.2]{Ath2}, which computes the face numbers of the Catalan arrangement. 
\begin{corollary} 
For $\lambda = (2, 2,\ldots,2)$ and $1 \leq k \leq d$, the number of lower faces of $\Pi_d (\lambda)$ of dimension $d-k$ is
\[ \sum_{r=k}^{d} \ (r-1)! \, S(d,r) \, {r \choose k} {r+k \choose k-1}, \]
where $S(d,k)$ stands for a Stirling number of the second kind.
\qed
\end{corollary}

\section{Remarks}

1. A direct bijective proof of Corollary 4.3 is possible and is outlined next. 

\vspace{0.1 in}
\noindent
\textit{Alternative proof of Corollary 4.3}: Recall that a region of ${\mathcal A}_d (\lambda)$ is defined by a linear ordering $\tau_1 < \tau_2 < \cdots < \tau_n$ of 
\[ x_i + m, \ \ 0 \leq m \leq \lambda_i - 1, \ 1 \leq i \leq d \]
which respects the orderings $x_i < x_i + 1 < \cdots < x_i + \lambda_i - 1$ for $1 \leq i \leq d$ and is such that the $\lambda$-permutation obtained by replacing each $x_i + m$ with $i$ is non-nesting. Given such an ordering $\tau$ with corresponding permutation $w = w_1 w_2 \cdots w_n$, let $j_1, j_2,\ldots,j_d$ be the positions in which $x_1, x_2,\ldots,x_d$ appear in $\tau$ respectively. In other words, $j_i$ is the smallest $j$ such that $w_j = i$. Consider the quotient of the abelian group ${\mathbb Z}_{n+1} ^d$ by the cyclic subgroup $H$ generated by $(1, 1,\ldots,1)$. The map
\[ \tau \toto J_{\tau} = (j_1, j_2,\ldots,j_d) + H \]
defines a bijection between the regions of ${\mathcal A}_d (\lambda)$ and cosets $J = (j_1, j_2,\ldots,j_d) + H$ of ${\mathbb Z}_{n+1} ^d / H$ for which all $j_i$ are mutually distinct. Clearly, the number of such cosets is $n (n-1) \cdots (n-d+2)$.
\qed

\begin{figure}[htpb]
\center{\mbox{\epsfbox{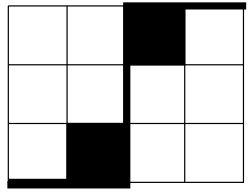}}} 
\caption{The subdivision defined by $(1, 1 2, 2, 1 2, 1)$}
\label{sub}
\end{figure}

\vspace{0.1 in}
2. Let $\pi : P \toto Q$ be as in the beginning of Section 4. The poset of $\pi$-coherent subdivisions of $Q$ is an induced subposet of the poset $\omega(P, Q)$ of all \textit{$\pi$-induced subdivisions} of $Q$ \cite[Definition 9.1]{Zi} (see also \cite[\S 3]{BKS}). In the special case $\dim Q = 1$ we are considering, this is the poset of cellular strings of $P$, induced by $\pi$ \cite[\S 1]{BKS}. The $\pi$-induced subdivisions in the lower part $P = \widetilde{\mathcal P}_{d+1} (\lambda)$ form an interesting poset, which we denote by $\Omega(\lambda)$. It is the set of all \textit{proper} ordered partitions of the multiset $M_{\lambda}$, partially ordered by refinement. An ordered partition of $M_{\lambda}$ is proper if none of its blocks contains repeated elements. The subdivision defined by $\rho = (B_1, B_2,\ldots,B_k) \in \Omega(\lambda)$ has as maximal faces the lifts of the faces
\[ \sum_{i=1}^{j-1} e_{B_i} + \sum_{r \in B_j} [0, e_r], \ \ 1 \leq j \leq k \]
of ${\mathcal P}_d (\lambda)$, where $e_B = \sum_{i \in B} e_i$. Figure \ref{sub} shows the faces of ${\mathcal P}_2 (4,3)$ which correspond to the proper ordered partition $(1, 1 2, 2, 1 2, 1)$ of $M_{\lambda} = \{1, 1, 1, 1, 2, 2, 2\}$, with $\lambda = (4, 3)$. The minimal elements of $\Omega(\lambda)$ are the $\lambda$-permutations, which give rise to the lifted monotone $\lambda$-paths.

\vspace{0.1 in}

The argument in the proof of Corollary 2.2 in Section 4 yields the following description of the subposet of $\pi$-coherent lower subdivisions: Let $\rho = (B_1, B_2,\ldots,B_k) \in \Omega(\lambda)$, as before. For each $1 \leq i \leq d$ replace the entries in the blocks of $\rho$ equal to $i$, from left to right, with $x_i, x_i +1,\ldots,x_i + \lambda_i - 1$ respectively, to get $\rho^{\prime} = (B_1^{\prime}, B_2^{\prime},\ldots,B_k^{\prime})$. The $\pi$-induced subdivision defined by $\rho$ is coherent if and only if the following linear system in the variables $x_1, x_2,\ldots,x_d$ is feasible:
\begin{enumit}
\item
$u = v$ whenever $u, v$ are in the same block $B_j^{\prime}$, 
\item
$u < v$ whenever $u \in B_j^{\prime}$, $v \in B_l^{\prime}$ and $j < l$.
\end{enumit}
In this case, the set of solutions forms exactly a face of ${\mathcal A}_d (\lambda)$. The subdivision of Figure \ref{sub} is not coherent because the system
\[ x_1 < x_1 + 1 = x_2 < x_2 + 1 < x_1 + 2 = x_ 2 + 2 < x_ 1 + 3 \]
has no solution. The only coherent atom of $\Omega(4,3)$ smaller than this subdivision is $(1, 2, 1, 2, 1, 2, 1)$.

\vspace{0.1 in}
3. There is a notion of ``flip'' on the set of all maximal monotone edge paths $p$ of $P$ with respect to the affine function $\pi$. Two such paths are said to be related by a \textit{flip} if they have a common cover in the poset $\omega(P, Q)$ of $\pi$-induced subdivisions of $Q$. We call the minimum number of flips required to make $p$ coherent the \textit{incoherency} of $p$. In the case of lifted monotone $\lambda$-paths $\gamma^w$, the operation of flipping swaps two distinct successive entries $i j$ of $w$. Therefore the incoherency of $\gamma^w$ is the minimum number of swappings needed to make $w$ non-nesting. To get a feeling for this number, we compute the maximum incoherency that can occur for $\lambda = (2, 2,\ldots,2)$.
\begin{proposition}
Let $\lambda = (2, 2,\ldots,2)$. The maximum incoherency of a path $\gamma^w$ is $d \choose 2$ and is attained by the $\lambda$-permutations of the form $w_{\sigma} = \sigma_1 \, \sigma_2 \cdots \sigma_d \, \sigma_d \cdots \sigma_2 \, \sigma_1$, where $\sigma = \sigma_1 \, \sigma_2 \cdots \sigma_d$ is a permutation of $\{1,\ldots,d\}$.
\end{proposition}
\begin{proof}
Every flip can reduce the number of ``nestings'' of $w$ by at most one and $w^{\sigma}$ has the maximum number $d \choose 2$ of nestings. On the other hand, any $w = w_1 w_2 \cdots w_{2d}$ can be turned into some $w_{\sigma}$ by at most $d \choose 2$ flips in the following way: Suppose that $w_1,\ldots,w_i$ are mutually distinct but $w_{i+1} = w_k$ for some $1 \leq k \leq i$. Let $j$ be the smallest index with $j > i$ for which $w_j \neq w_1,\ldots,w_i$. Move $w_j$ to the left until it occupies position $i+1$. Continue in the same way until the first $d$ entries of the permutation are all distinct, say $\sigma_1, \sigma_2,\ldots,\sigma_d$. Now continue flipping in the second half of the permutation, in an obvious way, to get $w_{\sigma}$ for $\sigma = \sigma_1 \, \sigma_2 \cdots \sigma_d$. Any two distinct integers $i, j$ are flipped at most once during the whole process, so the total number of flips is at most ${d \choose 2}$, as desired.  
\end{proof}

\vspace{0.1 in}
4. The extended Catalan arrangements form one of the two families of deformations of ${\mathcal A}_d$ for which explicit formulas have been obtained for the number of faces of any given dimension \cite[Ch.\ 8]{Ath2}. The other family consists of the \textit{extended Shi arrangements} 
\[ x_i - x_j = -a+1, -a+2,\ldots,a \ \ \for \ \ 1 \leq i < j \leq d, \]
where $a \geq 1$ is an integer. The number of faces in each dimension was shown to have a surprisingly simple combinatorial interpretation \cite[Thm.\ 8.2.1]{Ath2}. We don't know if a fiber polytope construction for the polar complexes exists in this case.

\vspace{0.1 in}
5. The cones of the extended Catalan arrangements ${\mathcal A}_d ^{[0,a]}$ were shown to be \textit{inductively free} \cite[Ch.\ 4]{OT} by Edelman and Reiner (see the proof of Theorem 3.2 in \cite{ER}). Other classes of deformations of ${\mathcal A}_d$, including the extended Shi arrangements, were shown to be inductively free in \cite{Ath3}. Proposition 4.2 suggests that the same is true for the arrangements ${\mathcal A}_d (\lambda)$.

\vspace{0.2 in}
\textit{Acknowledgement.} 
I am grateful to Lou Billera and Bernd Sturmfels for helpful discussions.

\vspace{0.1 in}

\end{document}